\documentclass[a4paper,11pt]{article}

\usepackage[english]{babel}
\usepackage[applemac]{inputenc}
\usepackage[OT1]{fontenc}

\usepackage{graphicx}
\usepackage[authoryear,round]{natbib}
\usepackage{amsmath,amsfonts,amssymb,amsthm}
\usepackage{mathrsfs}
\usepackage{verbatim}
\usepackage{numprint}
\usepackage{enumerate}

\newcommand{\Real}{\mathrm{Re}}
\newcommand{\prob}[1]{\operatorname{\mathbf{P}}\left(#1\right)}
\newcommand{\drv}{\mathrm{d}}
\newcommand{\Exp}{\mathrm{e}}
\newcommand{\expt}[1]{\operatorname{\mathbf{E}}\left[#1\right]}
\newcommand{\ind}[1]{\operatorname{\mathbf{1}}_{ \{ #1 \} } }
\newcommand{\eqdist}{\; {\buildrel d \over =} \;}
\newcommand{\wedgemin}[2]{#1 \wedge #2}

\DeclareMathOperator*{\argmax}{arg\,max}

\theoremstyle{plain}

\newtheorem{proposition}{Proposition}

\theoremstyle{remark}

\title{Calculating optimal limits for transacting credit card customers}
\author{J.\,K.~Budd \& P.\,G.~Taylor\thanks{School of Mathematics and
    Statistics, The University of Melbourne, Parkville 3010,
    Australia. Email:jkbudd@ms.unimelb.edu.au Email:taylorpg@unimelb.edu.au}}
\date{\vspace{-7ex}}

\begin{document}
\maketitle
\begin{abstract}
  We present a model of credit card profitability, assuming that the
  card-holder always pays the full outstanding balance. The motivation
  for the model is to calculate an optimal credit limit, which
  requires an expression for the expected outstanding balance. We
  derive its Laplace transform, assuming that purchases are made
  according to a marked point process and that there is a simplified
  balance control policy in place to prevent the credit limit being
  exceeded. We calculate optimal limits for a compound Poisson process
  example and show that the optimal limit scales with the distribution
  of the purchasing process and that the probability of exceeding the
  optimal limit remains constant. We establish a connection with the
  classic newsvendor model and use this to calculate bounds on the
  optimal limit for a more complicated balance control
  policy. Finally, we apply our model to real credit card purchase
  data.
\end{abstract}
\textbf{Keywords:} Banking; finance; optimization; stochastic
processes; statistics; time series. %
\section*{Introduction}
\label{sec:Introduction}
Managers of retail credit card portfolios regularly employ modelling
techniques to aid and automate decision-making throughout the customer
life-cycle. At the point where the customer is acquired, the decision
is whether or not to grant credit, and, if credit is to be granted,
what amount. Later in the customer life-cycle, incentives such as an
increased limit or a different interest rate may be offered, so the
problem is to determine those customers most likely to accept the
offer, or those who will generate the most profit if they accept. \par
There are numerous techniques available for modelling such
decisions. For a review of the most common techniques in use, see
\cite{rosenberg1994quantitative},
\cite{thomas2002credit,thomas2004readings}, \cite{crook2007recent}, or
\cite{hand1997statistical} for a particular emphasis on statistical
techniques. The particular problem of credit limit assignment has been
analysed by several authors. \cite{bierman1970credit} formulated a
dynamic programming model in which the decision variables were whether
or not to grant credit and what amount. In their formulation, the
amount of credit offered was linked to the probability of non-payment
via an exponentially declining relationship. Their model was later
extended in \cite{dirickx1976extension} to remove the assumption that
there is a zero expected future payoff from period $i$ onwards if no
payment is made. More recently, \cite{trench2003managing} developed a
Markov decision process (MDP) in which the objective was to optimise
customer lifetime value by either changing the customer's credit limit
or interest rate. \cite{so2011modelling} also used a MDP to generate a
dynamic credit limit policy where the state space was the account
behaviour score. \par
The problem that we address in this paper differs from those mentioned
above in that we are interested in the situation where credit has
already been granted and a limit assigned. A credit risk manager or
analyst may want to review a customer's limit if, for example, the
customer has contacted the card issuer to request an increased limit
or request a review of an automated decision to decline a limit
increase. We suppose that individual transaction data (the history of
purchases and payments) has been collected and retained and the
question is therefore whether we should revise an individual
customer's limit in light of their spending and payment behaviour with
the objective of maximising profitability for the bank. In particular,
we focus on the situation where the customer exhibits
\emph{transacting} behaviour. That is, the customer regularly pays the
outstanding balance in full by the due date issued on their monthly
statement. \par
We propose modelling at the individual level using transaction data as
a way to create further differentiation within certain populations
encountered in credit card portfolios. It is quite common to see large
numbers of low-risk accounts with similar behaviour scores and since
this is often a key input into account management strategies, these
accounts will receive identical treatment. Analysis of the
individual transaction patterns of such accounts can provide increased
differentiation and lead to the development of more profitable
strategies. \par
Indeed, the availability of transaction-level data has soared in the
last decade with the increasing prevalence of data warehouses in
financial institutions. Despite this increased availability, there
have been few attempts made to utilise this data to develop new models
for account management. \cite{hand2001prospecting} detailed the use of
data mining techniques to reveal spending patterns in a database of UK
credit card transactions. In particular, they modelled spending
behaviours at petrol stations. \cite{till2003behavioural} also used a
database of credit card transactions at petrol stations and fitted
various distributions to the inter-transaction times. This was
explored in more detail in \cite{tillThesis2001}. \par
The model developed in this paper bears similarities to those used in
inventory theory, in particular, the model analysed in
\cite{arrow1951optimal}, which is more commonly referred to as the
newsvendor model and originally attributed to
\cite{edgeworth1888mathematical}. The newsvendor model is known to
have a solution in terms of the quantile function of the input demand
distribution, but the differences introduced by our model formulation
require a solution by other means. We follow a method similar to that used in
\cite{chieraTaylor2002} to obtain a solution in terms of Laplace
transforms, which we then invert numerically following the algorithm
described in \cite{abate1995numerical}.
\section*{A model of individual credit card profitability}
\label{sec:creditCardProfitability}
Consider a credit card with limit $\ell > 0$ and no annual fees or
loyalty scheme. Let $B_\ell(t)$ denote the outstanding balance at
time $t$ and let $R(t)$ denote the cumulative profit earned up to time
$t$. We assume that the lending institution must pay a proportion of
the limit $\nu$ as a cost of financing. Now, let
$0 = t_0 < t_1 < \cdots < t_n$ be a sequence of billing times where
the outstanding balance is billed to the customer and
$s_1, \cdots, s_n$ be a sequence of times by which full or partial
payment of the outstanding balance is due, with each
$s_i \in (t_i , t_{i+1})$, $i = 1,\dots,n$. The times $s_i$ are
commonly referred to as due dates, the interval $(t_i, t_{i+1})$ as
the $i$th statement period, and the interval $(t_i,s_i)$ as the $i$th
interest-free period. To simplify matters, we assume that
$s_i = t_i, \, i \geq 1$, but note this does not affect the
generality of our model as we shall see.  \par
If we further assume that the customer exhibits transacting behaviour
and pays the full balance due before the end of each statement period,
then there will be no interest charges and the only contribution to
revenue will be from interchange which we assume occurs at a
proportion $\gamma$ of the total purchase amount. The only cost will
be the cost of funding the limit assigned to the customer and so the
profit earned in the period $(t_i, t_{i+1})$ is
\begin{equation}
  \label{eq:transactorProfit}
  R(t_{i+1}) - R(t_i) = \gamma B_\ell(t_{i+1}) - \nu \ell, \quad 0 \leq i
  \leq n,
\end{equation}
since the only contribution to the balance will be from new purchases in
the period $(t_i,t_{i+1})$. \par
Assume now that the each statement period has fixed length $T > 0$ and
that the customer's purchasing behaviour remains the same across each
period. Then we need only consider a single period and we may rewrite
\eqref{eq:transactorProfit} as
\begin{equation}
  \label{eq:transactorProfitFixedLength}
  R(T) = \gamma B_\ell(T) - \nu \ell.
\end{equation}
Now taking the expectation of \eqref{eq:transactorProfitFixedLength},
we have
\begin{equation}
  \label{eq:expectedTransactorProfit}
  \expt{R(T)} = \gamma \expt{B_\ell(T)} - \nu \ell.
\end{equation}
If there is an interest-free period of fixed length $b > 0$, then it is
not hard to see that we may account for its effect by adding this to the
length of the statement period. \par
We now seek to maximise the expected revenue by finding an appropriate
limit $\ell \in \Lambda$ where $\Lambda$ is a set of permitted
limits. Define
\begin{equation}
  \label{eq:maxCostFunction}
 \hat{\ell} := \argmax_{\ell \in \Lambda} \Big\{ \gamma \expt{B_\ell(T)} -
 \nu \ell \Big\}.
\end{equation}
If $\Lambda$ is a finite set, then we may determine $\hat{\ell}$
by simply evaluating the right-hand side of \eqref{eq:maxCostFunction} at
each point in $\Lambda$. The situation where  $\Lambda$ is not
countable requires some knowledge of the properties of
$\expt{B_\ell(T)}$. If $\expt{B_\ell(T)}$ happens to be a
differentiable function of $\ell$ at the point $T$, and if the maximum
in \eqref{eq:maxCostFunction} occurs in the interior of $\Lambda$, we can
determine $\hat{\ell}$ by differentiating
\eqref{eq:expectedTransactorProfit} and setting the right-hand side
equal to $0$. This yields
\begin{equation}
  \label{eq:transactorRevenueOptimisation1}
  \frac{\nu}{\gamma} = \frac{\partial}{\partial \ell} \expt{B_\ell(T)}.
\end{equation}
The task is now to find the limit $\hat{\ell}$ that will render the
right-hand side of \eqref{eq:transactorRevenueOptimisation1} equal to
the left-hand side. Whether we are solving \eqref{eq:maxCostFunction}
in the general case, or determining $\hat{\ell}$ via
\eqref{eq:transactorRevenueOptimisation1}, we require an expression
for $\expt{B_{\ell}(T)}$ or its derivative if it exists. We develop
such expressions in the next section.
\section*{An integral equation for the tail function}
\label{sec:integralEquation}
We assume that the card-holder attempts to make purchases according to
a marked point process
\begin{equation}
  \label{eq:mppDefinition}
  A(t) = \sum_{i = 1}^{N(t)} \xi_i,
\end{equation}
where $\{ \xi_i \}_{i=1,2,\dots}$ is a sequence of non-negative,
independent random variables with common distribution function $F$ and
$N(t)$ is a random variable, independent of $\{\xi_i\}_{i=1,2,\dots}$,
describing the number of events in $(0,t]$ in a renewal process with
inter-event time distribution $G$. For $k=1,2,\ldots$, we write
$t_k= \inf \{ t:N(t)=k \}$ and $\tau_k = t_k - t_{k-1}$, with
$t_0 = 0$. For the remainder of the paper, we assume that both $F$ and
$G$ are of exponential order and that all moments of the distributions
exist. These conditions are sufficient to ensure the existence of the
Laplace transforms
\begin{equation}
  \label{eq:distributionTransforms}
  \tilde{f}(\theta) = \int_0^\infty \! \Exp^{- \theta z} \, F(\drv z)
  \quad \mathrm{and} \quad \tilde{g}(\omega) = \int_0^\infty \! \Exp^{- \omega u} \, G(\drv u)
\end{equation}
for $\Real(\theta) > \sigma_F$ and $\Real(\omega) > \sigma_G$, where
the respective abscissae of convergence $\sigma_F$ and $\sigma_G$ of
$\tilde{f}$ and $\tilde{g}$ are strictly less than zero.
\par
Suppose now that the bank enforces a control policy on the outstanding
balance whereby if an attempted purchase would cause the outstanding
balance to exceed the credit limit $\ell$, that purchase is rejected
and the customer is barred from making any further purchases until the
outstanding balance is repaid at the end of the statement
period. Since our customer is a transactor, payment of the outstanding
balance in full is guaranteed. \par
We find an expression for the tail function $S_\ell(y,t) :=
\prob{B_\ell(t) \in (y,\ell] }$ by conditioning on the time and value
of the first jump of the process. For $0 < y \leq \ell$, we have three
possibilities to consider
\begin{enumerate}[i.]
\item The process jumps to some $z \in (0,y]$ at some time
  $u \in (0,t)$ and then regenerates itself at this point. That is, a
  new process starts at $z$ that behaves like the original one, but
  shifted by $\ell - z$ in space and $t - u$ in
  time. \label{enum:poss1}
\item The process jumps to some $z \in (y,\ell]$ and any subsequent
  jumps that happen in the remaining time interval $(u,t)$ cannot take
  the process out of the interval $(y,\ell]$. \label{enum:poss2}
\item The process jumps to some $z \in (\ell,\infty)$. If this occurs, the
  process is frozen at the point from which it jumped.
\end{enumerate}
With $\tau=\tau_1$ and $\xi=\xi_1$ we combine the cases above to
derive
\begin{equation}
  \label{eq:tailFunctionConditionalCases0}
  \expt{\ind{B_\ell(t) \in (y,\ell]} \mid \tau, \, \xi} =
  \begin{cases}
    \ind{B_{\ell}(t) - B_{\ell}(\tau) \in (y-\xi,\ell-\xi]}, & \tau \leq t, \,
     \xi \leq y \\
    1, & \tau \leq t, \, y < \xi \leq \ell
  \end{cases}
\end{equation}
By the regenerative property mentioned above, the distribution of
$B_{\ell}(t) - B_{\ell}(\tau)$ is the same as that of
$B_{\ell-\xi}(t-\tau)$, which assumes that the payment period is
$t-\tau$ and the credit limit is $\ell-\xi$. So we have
\begin{equation}
  \label{eq:tailFunctionConditionalCases}
\expt{\ind{B_\ell(t) \in (y,\ell]} \mid \tau, \, \xi} =
  \begin{cases}
    \ind{B_{\ell-\xi}(t-\tau) \in (y-\xi,\ell-\xi]}, & \tau \leq t, \,
    \xi \leq y \\
    1, & \tau \leq t, \, y < \xi \leq \ell
    \end{cases}
\end{equation}
By the law of total probability,
\begin{equation}
  S_\ell(y,t) = \int_0^t \int_0^y \! S_{\ell-z}(y-z,t-u) F(\drv z) \, G(\drv u)
  + G(t) \big( F(\ell) - F(y)
  \big). \label{eq:tailGenerator}
\end{equation}
We now wish to obtain the Laplace transform
\begin{equation}
 \tilde{S}(\theta,\omega,\psi) := \int_0^\infty \! \! \int_0^\infty \!
  \! \int_0^\ell \! \Exp^{-(\omega t + \theta \ell + \psi y)} \,
  S_\ell(y,t) \, \drv y \, \drv \ell \, \drv t. \label{eq:tailLST}
\end{equation}
It follows from \eqref{eq:tailGenerator} that
\begin{equation}
  \label{eq:exponentialOrderInequality}
  |S_\ell(y,t)| \leq G(t) F(y) + G(t) F(\ell) - G(t)F(y) = G(t) F(\ell),
\end{equation}
and since products of functions of exponential order are also of
exponential order, we have that $S_\ell(y,t)$ is of exponential order
when $F$ and $G$ are, and hence, the Laplace transform
\eqref{eq:tailLST} exists. \par
Applying the Laplace transform to \eqref{eq:tailGenerator},
we have after some rearrangement,
\begin{equation}
  \label{eq:tailGeneratorLST}
  \tilde{S}(\theta,\omega,\psi) = \frac{1}{\theta
      \omega \psi}   \frac{\tilde{g}(\omega) \big( \tilde{f}(\theta) -
      \tilde{f}(\theta + \psi) \big) }{1 -
      \tilde{g}(\omega)\tilde{f}(\theta + \psi)} .
\end{equation}
\par
To calculate the (two-dimensional) Laplace transform of the expectation, we note that
\begin{align}
  \label{eq:expectationOfNonNegRV}
  \mathcal{L}_{\theta,\omega} \big\{
  \expt{B_\ell(t)} \big\} & := \int_0^\infty \! \! \int_0^\infty \! \!
                            \Exp^{ - (\omega t + \theta \ell)}
                            \expt{B_\ell(t)} \drv \ell \, \drv t \\
                          & = \int_0^\infty \! \! \int_0^\infty \!
                            \Exp^{ - (\omega t + \theta \ell)}
                            \! \int_0^\ell \! \prob{B_\ell(t) \in (y,\ell]} \drv y \, \drv
                            \ell \, \drv t,
\end{align}
which corresponds to evaluating $\tilde{S}(\theta,\omega,0)$. We
apply l'H\^opital's rule to \eqref{eq:tailGeneratorLST} to obtain
\begin{align}
  \lim_{\psi \rightarrow 0}
  \tilde{S}(\theta,\omega,\psi) & =
                                  \lim_{\psi \rightarrow 0}
                                  \frac{1}{\theta
                                  \omega
                                  \psi}
                                  \frac{\tilde{g}(\omega)
                                  \big(\tilde{f}(\theta)
                                  -
                                  \tilde{f}(\theta
                                  +
                                  \psi)\big)}{1
                                  -
                                  \tilde{g}(\omega)\tilde{f}(\theta
                                  + \psi)}
                                  \nonumber
  \\
                                & = -
                                  \frac{\tilde{g}(\omega)}
                                  {\theta \omega \big( 1
                                  - \tilde{g}(\omega)
                                  \tilde{f}(\theta)
                                  \big)} \frac{\drv}{\drv
                                  \theta} \,
                                  \tilde{f}(\theta) \label{eq:expectationLST}.
\end{align}
It should be noted that the derivative of $\expt{B_\ell(t)}$ may not
exist. Indeed, if $F$ is lattice then $\expt{B_\ell(t)}$ will be a
step function. In the case where the derivative does exist, we obtain
its Laplace transform by multiplying \eqref{eq:expectationLST} by
$\theta$ to yield
\begin{equation}
  \mathcal{L}_{\theta,\omega} \bigg\{
  \frac{\partial}{\partial \ell} \expt{B_\ell(t)}
  \bigg\} = - \frac{\tilde{g}(\omega)}{\omega \big(1 -
    \tilde{g}(\omega)\tilde{f}(\theta) \big) } \frac{\drv}{\drv
    \theta} \tilde{f}(\theta). \label{eq:derivativeLST}
\end{equation}
\section*{An example using a compound Poisson process}
\label{sec:example}
\addcontentsline{toc}{section}{An example using a compound Poisson
  process} In this section we use Equation \eqref{eq:derivativeLST} to
calculate the optimal limit for a transacting credit card customer who
makes purchases according to a compound Poisson process with rate
$\lambda$ and purchase sizes that are exponentially distributed with
parameter $\mu$. The Laplace transforms of the inter-arrival time
distribution and the purchase size distribution are
\begin{equation}
 \tilde{f}(\theta) = \frac{\mu}{\mu + \theta}, \quad
 \mathrm{and} \quad \tilde{g}(\omega) = \frac{\lambda}{\lambda + \omega}
\end{equation}
so now equation \eqref{eq:derivativeLST} becomes
\begin{equation}
   \mathcal{L}_{\theta,\omega} \bigg\{ \frac{\partial}{\partial \ell} \expt{B_\ell(t)}
  \bigg\} = \frac{\lambda \mu}{\omega (\theta + \mu) \big( \mu
    \omega + \theta (\lambda + \omega) \big)}, \label{eq:cppDerivativeLST}
\end{equation}
which can be inverted analytically to yield
\begin{equation}
  \mathcal{L}_{\theta} \bigg\{ \frac{\partial}{\partial \ell} \expt{B_\ell(t)}
  \bigg\} = \frac{1}{\theta} \bigg( \frac{\mu}{\mu + \theta} \bigg) \Bigg( 1 -
  \exp \bigg\{ \lambda t \bigg(
  \frac{\mu}{\mu + \theta} - 1 \bigg) \bigg\} \Bigg). \label{eq:cppDerivativeLSTInverted}
\end{equation}
It does not appear to be easy to perform further analytical inversion
with respect to $\theta$ since the inverse transform of the
exponential of a rational function is not a standard transform. As
such, we resort to numerical inversion using the \texttt{EULER}
algorithm as detailed in \cite{abate1995numerical}. \par
We calculated the optimal limit using an interchange rate
$\gamma = 0.0054$, a cost of funds $\nu = 0.0007$ and a statement
period length $T = 30$. We took $\Lambda = (0,5000]$ and used a
bisection search to solve Equation
\eqref{eq:transactorRevenueOptimisation1}. Table
\ref{tab:optimalLimitResults} shows the optimal limits calculated for
a range of values of $\lambda$ and $\mu$. \par
An important factor in the assignment of credit limits is the
customer's experience of having a purchase declined due to
insufficient funds which, in our model, will result in the customer
being barred from making further purchases until the end of the
statement period. The probability of a purchase being declined is
given by the tail function of $A(T)$ which, in the case of a compound
Poisson process with exponential marks, has Laplace transform
\begin{equation}
  \label{eq:cppTailLST}
  \tilde{A}(\psi) = \frac{1}{\psi} \Bigg( 1 - \exp \bigg\{ \lambda T \Big(
    \frac{\mu}{\mu + \psi} - 1 \Big) \bigg\} \Bigg), \quad \Real(\psi) > -\mu.
\end{equation}
The results of inverting Equation \eqref{eq:cppTailLST} at the optimal
limits calculated in Table \ref{tab:optimalLimitResults} are presented
in Table \ref{tab:optimalLimitBlockingProb}. They show that the
probability of a declined purchase remains constant when the rate
parameter of the purchase size distribution changes. Indeed, by
scaling the purchase size distribution by some
$\alpha \in \mathbb{R}_+$, we scale the input process $A(t)$ and, as
evidenced by Table \ref{tab:optimalLimitResults}, the optimal
limit. The following proposition formalises this result.
\begin{proposition}[Scaling property of the optimal limit]
  \label{prop:scaling}
  Let $A'(t) \eqdist \alpha A(t)$ and let
  \begin{equation*}
    B_\ell(t) = \sup_{0 \leq u \leq t} \{A(u) : A(u) \leq \ell \} \quad
    \mathrm{and} \quad
    B_{\ell}'(t) = \sup_{0 \leq u \leq t} \{A'(u) : A'(u) \leq \ell \}.
  \end{equation*}
Then the solution to the optimisation problem
\begin{equation*}
 \hat{\ell}' =  \argmax_{\ell \in \Lambda} \big( \gamma \expt{B_\ell'(t)} - \nu \ell
  \big) = \alpha \argmax_{\ell \in \Lambda} \big( \gamma \expt{B_\ell(t)} - \nu \ell
  \big) = \alpha \hat{\ell},
\end{equation*}
as long as $\alpha \hat{\ell} \in \Lambda$. Furthermore, we have that
\begin{equation*}
  \prob{A'(t) > \hat{\ell}'} = \prob{A(t) > \hat{\ell}}.
\end{equation*}
\end{proposition}
The proof of the above proposition is given in the appendix.
\npdecimalsign{.}
\nprounddigits{2}
\begin{table}[htbp]
  \centering
  \begin{tabular}{cc|rrrrr}
    & & \multicolumn{5}{c}{$1/\mu$} \\
    & & 20 & 40 & 60 & 80 & 100 \\
    \hline
    & 1 & \numprint{798.2182713702} & \numprint{1596.4365427400} & \numprint{2394.6548141093} & \numprint{3192.8730854797} & \numprint{3991.0913568496} \\
    & 2 & \numprint{1470.4018320563} & \numprint{2940.8036641135} & \numprint{4411.2054961706} & \numprint{5881.6073282277} & \numprint{7352.0091602848} \\
    $\lambda$ & 3 & \numprint{2125.8575397851} & \numprint{4251.7150795709} & \numprint{6377.5726193557} & \numprint{8503.4301591393} & \numprint{10629.2876989269} \\
    & 4 & \numprint{2772.6317142248} & \numprint{5545.2634284508} & \numprint{8317.8951426747} &\numprint{ 11090.5268568997} & \numprint{13863.1585711252} \\
    & 5 & \numprint{3413.8515798266} & \numprint{6827.7031596528} & \numprint{10241.5547394791} & \numprint{13655.4063193053} & \numprint{17069.2578991332}
  \end{tabular}
  \caption{\small Table of values for the optimal limit. The values
    were calculated using a statement period length of $T=30$, with
    $\gamma=0.0054$ and $\nu = 0.0007$.}
  \label{tab:optimalLimitResults}
\end{table}
\nprounddigits{8}
\begin{table}[htbp]
  \centering
  \begin{tabular}{cc|rrrrr}
    & & \multicolumn{5}{c}{$1/\mu$} \\
    & & 20 & 40 & 60 & 80 & 100 \\
    \hline
    & 1 & \numprint{0.1059165813} & \numprint{0.1059165813} & \numprint{0.1059165813} & \numprint{0.1059165813} & \numprint{0.1059165813} \\
    & 2 & \numprint{0.1121867114} & \numprint{0.1121867114} &
                                                            \numprint{0.1121867114} & \numprint{0.1121867114} & \numprint{0.1121867114} \\
    $\lambda$ & 3 & \numprint{0.1151312748} & \numprint{0.1151312748} & \numprint{0.1151312748} & \numprint{0.1151312748} & \numprint{0.1151312748} \\
    & 4 & \numprint{0.1169382056} & \numprint{0.1169382056} & \numprint{0.1169382056} & \numprint{0.1169382056} & \numprint{0.1169382056} \\
    & 5 & \numprint{0.1181941315} & \numprint{0.1181941315} &
                                                            \numprint{0.1181941315}
           & \numprint{0.1181941315} & \numprint{0.1181941315}
  \end{tabular}
  \caption{\small Probability of the credit card customer experiencing
    a declined purchase when assigned the optimal limit.}
  \label{tab:optimalLimitBlockingProb}
\end{table}
\section*{Comparison with the newsvendor model}
\label{sec:newsvendor}
As mentioned in the introduction, the model we have formulated is
similar to the single period newsvendor model with random demand. The
objective of the newsvendor problem is to determine the number of
newspapers to stock which will maximise the expected profit. Following
the formulation in \cite{porteus2002foundations}, let $A(T)$ denote
the random demand for newspapers in a single period of length $T$,
$\ell$ the stock level of newspapers and $\gamma$ and $\nu$ the unit
profit and cost respectively. The expected revenue earned by the
newsvendor in a period is $\gamma \expt{\wedgemin{A(T)}{\ell}}$,
representing the cases where the newsvendor has either ordered
sufficiently-many newspapers to meet the demand $A(T)$, or an
insufficient number, in which case he sells the entire stock $\ell$
ordered at the beginning of the period. The cost incurred by the
newsvendor is simply the unit cost of each newspaper, $\nu$,
multiplied by the number of newspapers $\ell$ which he chooses to
order. Thus, the problem is to determine
\begin{equation}
  \label{eq:newsvendorProblem}
  \ell^* := \argmax_{\ell \in \mathbb{R}^+} \gamma
  \expt{\wedgemin{A(T)}{\ell}} - \nu \ell,
\end{equation}
which is similar to the problem formulated in Equation
\eqref{eq:maxCostFunction}. The newsvendor problem has a well-known
solution in terms of the distribution of $A(T)$. For comparison, we
can rewrite \eqref{eq:maxCostFunction} as
\begin{equation}
  \label{eq:controlSystemOptimality3}
  \hat{\ell} = \argmax_{\ell \in \Lambda} \gamma \expt{ \wedgemin{A(T)}{\ell - U} } - \nu \ell, \\
\end{equation}
where $U$ is a random variable describing the undershoot
of the process, conditional on the event that $A(T) > \ell$. Since
an undershoot only occurs when the input process $A(t)$ exceeds $\ell$
before the end of the period $T$, we have
\begin{align}
  \expt{\wedgemin{A(T)}{\ell - U}} & = \prob{A(T) \leq \ell} \expt{A(T)
                                     \mid A(T) \leq \ell} \nonumber \\
                                   & \quad + \prob{A(T) > \ell}
                                     \expt{\ell - U \mid A(T) >
                                     \ell}. \label{eq:splitExpectation}
\end{align}
However, the solution using this formulation requires explicit
knowledge of the distribution of $U$ (or equivalently, $\ell-U$)
which is, in general, difficult to obtain. \par
The difference between the two models lies in the fact that, in
the case where a customer exceeds his or her credit limit, the balance
of purchases made is (with probability one) strictly less than the
limit $\ell$, whereas the newsvendor always sells $\ell$ papers
whenever the demand exceeds $\ell$. The latter model would apply to
the credit limit case if, whenever a customer attempts a purchase of
value $z$ that will take the current outstanding balance $x$ over the
credit limit $\ell$, only $\ell - x$ is charged to the card and then
no further purchases are allowed. This is clearly an unsatisfactory
rule to use in our model since this would involve merchants only
accepting partial payment for whatever goods were being
purchased. \par
An improved model for credit card use would allow for the customer to
continue to attempt purchases after a purchase has been declined since,
in reality, a bank's card management system will permit any number of
purchases to be made so long as the total value does not exceed
$\ell$.
Let $\bar{B}_\ell(t)$ denote the value of the outstanding balance
under this control policy. We can derive an integral equation for the
tail function of $\bar{B}_\ell(t)$ by adding another case to Equation
\eqref{eq:tailFunctionConditionalCases0} to include the possibility
that the process restarts from its original position if the first jump
takes the process over $\ell$. However, a closed-form expression for
the Laplace transform of the tail function of $\bar{B}_\ell(t)$ is not
immediately forthcoming when we include this third term. \par
Alternatively, we can obtain bounds on the limit that would be set
when further purchases are allowed following rejection by using the
newsvendor model and the model we have developed which prevents
further purchases after the first declined purchase. We claim
\begin{equation}
  \expt{B_\ell(T)} \leq \expt{\bar{B}_\ell(T)} \leq
  \expt{\wedgemin{A(T)}{\ell}}. \label{eq:controlSystemInequality}
\end{equation}
That $\expt{B_\ell(T)} \leq \expt{\wedgemin{A(T)}{\ell}}$ follows from
a direct comparison of the expectations in their integral form,
\begin{equation}
  \int_0^\ell z \, F_A(\drv z) + \big( 1-F_A(\ell) \big) \int_0^\ell z
  \, F_Z(\drv z) \leq \int_0^\ell z \, F_A(\drv z) + \big( 1-F_A(\ell) \big)
  \ell.
\end{equation}
We reason that $\expt{B_\ell(T)} \leq \expt{\bar{B}_\ell(T)}$ since a
balance control policy that allows for further purchases following a
rejected purchase cannot decrease the expected balance. Similarly,
$\expt{\bar{B}_\ell(T)} \leq \expt{\wedgemin{A(T)}{\ell}}$ since under
the newsvendor control policy, a rejected purchase will always result
in an outstanding balance of $\ell$, but this is not necessarily so
under the policy allowing for purchase retrials. \par
We now claim
\begin{equation}
  \label{eq:1}
  \hat{\ell} \geq \bar{\ell} \geq \ell^*
\end{equation}
where
\begin{align}
  \bar{\ell} & := \argmax_{\ell \in \Lambda} \big\{ \gamma \expt{\bar{B}_\ell(T)} - \nu \ell
               \big\} \label{eq:retrialProgram}
  \intertext{and}
               \ell^* & := \argmax_{\ell \in \Lambda} \big\{ \gamma
                        \expt{\wedgemin{A(T)}{\ell}} - \nu \ell
                        \big\}, \label{eq:newsvendorProgram}
\end{align}
which follows directly from \eqref{eq:controlSystemInequality}. \par
The optimal limit in the newsvendor model is given by
\begin{equation}
  \label{eq:newsvendorSolution}
  \ell^* = \inf \Big\{ \ell : F_A(\ell) \geq \frac{\gamma - \nu}{\gamma} \Big\},
\end{equation}
and Table \ref{tab:newsvendorLimitResults} shows the optimal limits
calculated when the input process $A(t)$ is
a compound Poisson process with arrival rate $\lambda$ and exponential
jumps with parameter $\mu$. The differences between the
limits obtained using the newsvendor model and the process $B_\ell(t)$
are shown in Table \ref{tab:newsvendorLimitDifferences}. These
differences give a measure of the model error due to using the balance
control policy where further purchases are prevented following an
attempt to exceed the credit card limit.
\begin{table}[htbp]
  \centering
  \begin{tabular}{cc|ccccc}
    & & \multicolumn{5}{c}{1/$\mu$} \\
    & & $20$ & $40$ & $60$ & $80$ & $100$ \\
    \hline
    &  $1$ & $776.78$ & $1,553.55$ & $2,330.33$ & $3,107.10$ & $3,883.88$ \\
    &  $2$ & $1,449.38$ & $2,898.76$ & $4,348.14$ & $5,797.52$ & $7,246.90$ \\
    $\lambda$ &  $3$ & $2,105.02$ & $4,210.04$ & $6,315.06$ & $8,420.09$ & $10,525.11$ \\
    &  $4$ & $2,751.91$ & $5,503.81$ & $8,255.72$ & $11,007.63$ & $13,759.54$ \\
    &  $5$ & $3,393.20$ & $6,786.41$ & $10,179.61$ & $13,572.81$ & $16,966.02$
  \end{tabular}
  \caption{\small Table of values for $\ell^*$, the optimal limit found
    using the newsvendor model, with $\gamma = 0.0054$, $\nu = 0.0007$
  and $T=30$.}
  \label{tab:newsvendorLimitResults}
  \begin{tabular}{cc|ccccc}
    & & \multicolumn{5}{c}{1/$\mu$} \\
    & & $20$ & $40$ & $60$ & $80$ & $100$ \\
    \hline
    & $1$ & $21.44$ & $42.89$ & $64.33$ & $85.77$ & $107.21$ \\
    & $2$ & $21.02$ & $42.05$ & $63.07$ & $84.09$ & $105.11$ \\
    $\lambda$  & $3$ & $20.84$ & $41.67$ & $62.51$ & $83.34$ & $104.18$ \\
    & $4$ & $20.72$ & $41.45$ & $62.17$ & $82.90$ & $103.62$ \\
    & $5$ & $20.65$ & $41.30$ & $61.95$ & $82.59$ & $103.24$
  \end{tabular}
  \caption{\small Table of values for $\hat{\ell} - \ell^*$, the
    difference between the optimal limits calculated using the process
    $B_\ell(t)$ and the newsvendor model.}
  \label{tab:newsvendorLimitDifferences}
\end{table}
\section*{An example using credit card transaction data}
\label{sec:real-data}
In this section, we apply the model we have developed to actual data
from a credit card customer. Two datasets of anonymised credit card
transactions were made available to the authors for the purposes of
this research. The first dataset holds posted transactions, which are
the approved purchases and payments, and also includes interest
charges, fees, reversals and other automated transactions. The second
dataset describes authorisations, which are the purchases and payments
attempted by customers. \par
The posted transactions dataset describes the value and processing
dates of $771,457$ transactions made between 8 February 2011 and 27
February 2013 by $3,734$ customers holding $3,971$ accounts. Of the
$771,457$ transactions, $511,969$ are retail purchase transactions and
$84,503$ are payments. In addition to the above, the dataset also
contains identifying merchant information which allows us to
categorise transactions by store type. \par
The authorisations dataset describes the value and transaction times,
accurate to the second, of $405,844$ transactions made between 7
February 2011 and 27 February 2013. The dataset also contains account
credit limits and describes whether or not the transaction was
approved or declined and, in the case of a decline, a code describing
the reason (e.g.\ insufficient funds or an incorrectly entered
PIN). These transactions were made by the same $3,734$ customers
across $4,333$ credit card accounts. Due to issues
encountered during the data extraction process, we were only able to
match the authorisation transaction records with the posted
transaction records of $2,246$ customers. These customers attempted
$288,423$ purchases, of which $223,804$ were approved. \par
For the purposes of illustrating the model developed in this paper, we
extracted the transactions of a single customer who was identified as
a transactor through the absence of interest charges to their account
over the period. We filtered the transactions to include only those
made at supermarkets since they account for a large proportion of
purchases made on the credit card ($306$ out of $732$) and are easily
identified in both the authorisations and posted transactions
dataset. The time series was modified to exclude transactions that
were declined due to a POS device error or an incorrect PIN entry. A
preliminary analysis of the supermarket transactions of several
card-holders revealed occasional clustering of transactions in
time. This could be explained by a number of customer behaviours. For
example, a customer may visit a supermarket only to find that some of
the items they intended on purchasing are not available, so they buy
the items that are in stock and then visit another supermarket nearby
to purchase the remaining items. It could also be due to a customer
forgetting some items, and quickly returning to the same store to
purchase them. With this in mind, transactions made within an hour of
each other were combined into a single transaction with the total
value of those transactions. \par
The customer made $306$ purchases at various supermarkets over a
period of $473$ days which totalled $\$11,469.44$. This equates to
approximately $\$37.36$ per transaction or $\$24.25$ per day. In a
$30$-day period, this totals $\$727.50$ in purchases, which is far
less than the account credit limit of $\$5,000$. \par
We fitted a $\Gamma$-distribution to the purchase values of the
modified time series and estimated the shape and scale parameters
using maximum likelihood estimation. Using the two-sided
Kolmogorov-Smirnov test statistic
\begin{equation}
  \label{eq:kstest}
  D_n = \sup_x | F_n(x) - F(x) |
\end{equation}
where $F_n(x)$ is the empirical distribution function and $F(x)$ is
the distribution function of the fitted $\Gamma$-distribution, we
found the fit to be statistically significant at the $0.05$ level as
evidenced by the result in Table \ref{tab:Stats}. Finding an
appropriate distribution for the inter-transaction times was not so
straight-forward, so for the purposes of this example, we assume the
inter-transaction times follow an exponential distribution with
parameter $\lambda$, which we estimated from the reciprocal of the
mean of the inter-purchase times to be
$\hat{\lambda} = 0.6451 \, \pm 0.0369$. We further assumed
independence of the purchase values and the inter-purchase times, but
note that this assumption could be tested by computing the coherence
between the inter-purchase times and the purchase values (see theorem
4.4 in \cite{brillinger2012spectral}). Some degree of dependence
between inter-purchase time and purchase value is likely, particularly
with supermarket transactions, since a large inter-purchase time would
indicate that a customer has not visited a supermarket for a while,
and hence the next purchase is likely to be a large one.
\begin{table}[htbp]
  \centering
  \begin{tabular}[c]{lr}
    \hline
    Statistic & Estimate \\
    \hline
    $D_n$ & $0.0350$ \\
    $p$-value & $0.8623$ \\ \hline
    $\hat{\mu}$ (shape) & $2.8946 \, \pm 0.2258$ \\
    $\hat{k}$ (scale) & $0.0769 \, \pm 0.0065$ \\
    \hline
  \end{tabular}
  \caption{Kolmogorov-Smirnov test statistics and
    $\Gamma$-distribution shape and scale parameter estimates for the
    purchase value distribution.}
  \label{tab:Stats}
\end{table}
Substituting
\begin{equation}
  \label{eq:6}
 \tilde{g}(\omega;\lambda) = \Bigg(
 \frac{\lambda}{\lambda + \omega} \Bigg) \quad \mathrm{and}
 \quad \tilde{f}(\theta;\mu,k) = \Bigg( \frac{\mu}{\mu +
    \theta} \Bigg)^{k}
\end{equation}
into Equation \eqref{eq:derivativeLST} and inverting once from
$\omega$ to $t$, we have for the Laplace transform of the expectation
and its derivative
\begin{equation}
  \label{eq:gammaExptLST}
 \mathcal{L}_{\theta} \big\{ \expt{B_\ell(t)} \big\} = \frac{k}{\theta(\mu + \theta)} \Bigg(
 \frac{\mu}{\mu + \theta}
 \Bigg)^k \frac{1 - \Exp^{\lambda t
       \left(\left(\frac{\mu }{\mu + \theta}\right)^k - 1
   \right)}}{1 - \left( \frac{\mu }{\mu + \theta} \right)^k}
\end{equation}
and
\begin{equation}
  \label{eq:gammaDrvLST}
 \mathcal{L}_{\theta} \bigg\{ \frac{\partial}{\partial \ell}
 \expt{B_\ell(t)} \bigg\} = \frac{k}{\mu + \theta} \Bigg(
 \frac{\mu}{\mu + \theta}
 \Bigg)^k \frac{1 - \Exp^{\lambda t
       \left(\left(\frac{\mu }{\mu + \theta}\right)^k - 1
   \right)}}{1 - \left( \frac{\mu }{\mu + \theta} \right)^k}.
\end{equation}
For calculations using the newsvendor model, we use the
Laplace transform of the tail function of the compound Poisson process
with $\Gamma$-distributed jumps,
\begin{equation}
  \label{eq:GammaCPPLST}
  \tilde{A}_\Gamma(\psi) =  \frac{1}{\psi} \Bigg( 1 - \exp \Bigg\{ \lambda T \left(
   \left( \frac{\mu}{\mu + \psi} \right)^k - 1 \right) \Bigg\} \Bigg),
 \quad \Real(\psi) > -\mu
\end{equation}
and
\begin{equation}
  \label{eq:minExpLST}
  \mathcal{L}_\theta \big\{ \expt{ \wedgemin{A(T)}{\ell} } \big\} =
  \int_0^\infty \! \Exp^{-\theta \ell} \! \int_0^\ell \prob{A(T) > y}
  \drv y \, \drv \ell  = \frac{1}{\theta} \tilde{A}_\Gamma(\theta).
\end{equation}
Again, we assume an interchange rate $\gamma = 0.0054$, cost of funds
$\nu = 0.0007$ and statement period length $T = 30$. Substituting the
estimated parameters $\hat{\lambda}$, $\hat{\mu}$ and $\hat{k}$ into
Equations \eqref{eq:gammaDrvLST}--\eqref{eq:GammaCPPLST}, we obtain
the results in Table \ref{tab:calibrationResults}. We again used a
bisection search and the \texttt{EULER} algorithm to calculate the
optimal limits. The table shows the expected balance, expected profit
and probability of a declined purchase at the original limit, the
upper and lower bounds of the optimal limit and a revised limit. The
bounds on the optimal limit accord with the average monthly
supermarket spend of the customer. Recall that we stated that the
customer spent $\$727.50$ in a $30$-day period; the upper and lower
limits yield an expected balance of just over $\$714$. The increase in
profitability is substantial, but we note that this is somewhat
artificial given we have restricted our analysis to only those
purchases made at supermarkets. The revised limit is proposed since
most card-issuers offer limits in multiples of $\$500$. As shown in
the table, the deviation from profit at optimality is negligible, but
there is slightly smaller chance of the customer experiencing a
declined purchase.
\begin{table}[htbp]
  \centering
  \begin{tabular}[h]{l|rrr}
    & Original & Optimal & Revised \\
    \hline
    Limit  & $\$5,000.00$ & $[\$947.83,\$973.81]$  & $\$1,000.00$ \\
    Expected balance & $ \$ 728.64$ & $[ \$ 714.06,\$ 714.09]$ & $[ \$ 717.13,
    \$ 719.64]$ \\
    Expected profit & $ \$ 0.44$ & $[ \$ 3.17, \$3.19] $ & $[ \$ 3.17,
                                                           \$3.19]$ \\
    Probability of decline & $0.0000$ & $[0.1060,0.1296]$ & $0.0857$
  \end{tabular}
  \caption{Expected balance, expected profit and probability of a declined
    purchase at the original limit, upper and lower bounds of the optimal
    limit and a proposed revised limit.}
  \label{tab:calibrationResults}
\end{table}
Although the results in Table \ref{tab:calibrationResults} show a
marked increase in profitability as a result of lowering the credit
limit to the optimal limit, it should be noted that doing so would
substantially increase the probability that the customer will
experience a declined purchase if their purchasing behaviour remains
unchanged. This is undoubtedly a poor experience for the customer and
the consequences of this for both the customer and the bank should be
considered before any change to the customer's credit limit is made.
\section*{Discussion}
\label{sec:discussion}
The model we have presented makes a number of simplifying
assumptions. As mentioned in the introduction, the assumption of transacting
behaviour is both valid and useful since most credit card portfolios
are primarily composed of customers exhibiting this behaviour and they
form a significant source of revenue through interchange. We can
extend the model to include the possibility of partial repayment of
the outstanding balance by including another term in
\eqref{eq:transactorProfit} which describes how partial repayment
generates interest. To then derive the resulting optimal limit
requires additional assumptions on payment behaviour and the value of
new purchases when the account retains a partially unpaid balance. We
regard this as a worthwhile avenue for future research given its
applicability in credit management. \par
We assumed that the attempted purchase process was a marked point
process with inter-purchase time distribution $G$ and purchase size
distribution $F$, and that these distributions remain unchanged with
respect to the outstanding balance and the credit limit. A more
realistic model would include state-dependence, as our data indicates
that customers either reduce their purchase frequency or purchase size
as they near their credit limit. Another modification would be to
split $A(t)$ into a series of marked point processes, each modelling
different types of purchases such as retail, restaurants,
supermarkets, or cash advances. \par
Related to the above, our calculation of the optimal limit assumes
that $A(t)$ will remain unchanged in the event of a change in the
credit limit: an improved model may factor in the reward associated
with an increased limit (customers may increase their overall
attempted spend) or, correspondingly, a penalty associated with a
decrease (a customer may decide to cancel their card). It was
demonstrated in \cite{soman2002creditlimit} that increasing credit
limits resulted in increased customer spend only in some portfolio
segments. However, there is relatively little published research on
the effect of credit limit decreases on customer purchasing
behaviour. Also, implementing the resulting optimal limit exactly may prove
difficult as most banks only offer limits to customers that are
multiples of $\$500$, for both customer experience and systems
reasons. Nonetheless, we regard the model as a useful complement to
existing limit setting strategies for understanding the effect of
limit changes on profitability and customer experience.
\section*{Acknowledgments}
The authors are grateful to Peter Braunsteins, Shaun McKinlay and
Nicholas Read for useful and enthusiastic discussions throughout the
research for this paper. P.\,G.~Taylor's research is supported by the
Australian Research Council (ARC) Laureate Fellowship FL130100039 and
the ARC Centre of Excellence for Mathematical and Statistical
Frontiers (ACEMS).
\bibliographystyle{apalike}
\bibliography{transactor_optimal_limits}

\begin{thebibliography}{}

\bibitem[Abate and Whitt, 1995]{abate1995numerical}
Abate, J. and Whitt, W. (1995).
\newblock Numerical inversion of {Laplace} transforms of probability
  distributions.
\newblock {\em ORSA Journal on Computing}, 7(1):36--43.

\bibitem[Arrow et~al., 1951]{arrow1951optimal}
Arrow, K.~J., Harris, T., and Marschak, J. (1951).
\newblock Optimal inventory policy.
\newblock {\em Econometrica: Journal of the Econometric Society}, pages
  250--272.

\bibitem[Bierman and Hausman, 1970]{bierman1970credit}
Bierman, Jr, H. and Hausman, W.~H. (1970).
\newblock The credit granting decision.
\newblock {\em Management Science}, 16(8):B--519.

\bibitem[Brillinger, 2012]{brillinger2012spectral}
Brillinger, D.~R. (2012).
\newblock The spectral analysis of stationary interval functions.
\newblock In Guttorp, P. and Brillinger, D.~R., editors, {\em Selected Works of
  David Brillinger}, pages 25--55. Springer.

\bibitem[Chiera and Taylor, 2002]{chieraTaylor2002}
Chiera, B. and Taylor, P. (2002).
\newblock What is a unit of capacity worth?
\newblock {\em Probability in the Engineering and Informational Sciences},
  16(04):513--522.

\bibitem[Crook et~al., 2007]{crook2007recent}
Crook, J.~N., Edelman, D.~B., and Thomas, L.~C. (2007).
\newblock Recent developments in consumer credit risk assessment.
\newblock {\em European Journal of Operational Research}, 183(3):1447--1465.

\bibitem[Dirickx and Wakeman, 1976]{dirickx1976extension}
Dirickx, Y.~M. and Wakeman, L. (1976).
\newblock An extension of the {Bierman-Hausman} model for credit granting.
\newblock {\em Management Science}, 22(11):1229--1237.

\bibitem[Edgeworth, 1888]{edgeworth1888mathematical}
Edgeworth, F.~Y. (1888).
\newblock The mathematical theory of banking.
\newblock {\em Journal of the Royal Statistical Society}, 51(1):113--127.

\bibitem[Hand and Blunt, 2001]{hand2001prospecting}
Hand, D.~J. and Blunt, G. (2001).
\newblock Prospecting for gems in credit card data.
\newblock {\em IMA Journal of Management Mathematics}, 12(2):173--200.

\bibitem[Hand and Henley, 1997]{hand1997statistical}
Hand, D.~J. and Henley, W.~E. (1997).
\newblock Statistical classification methods in consumer credit scoring: a
  review.
\newblock {\em Journal of the Royal Statistical Society: Series A (Statistics
  in Society)}, 160(3):523--541.

\bibitem[Porteus, 2002]{porteus2002foundations}
Porteus, E.~L. (2002).
\newblock {\em Foundations of stochastic inventory theory}.
\newblock Stanford University Press.

\bibitem[Rosenberg and Gleit, 1994]{rosenberg1994quantitative}
Rosenberg, E. and Gleit, A. (1994).
\newblock Quantitative methods in credit management: a survey.
\newblock {\em Operations research}, 42(4):589--613.

\bibitem[So and Thomas, 2011]{so2011modelling}
So, M. and Thomas, L.~C. (2011).
\newblock Modelling the profitability of credit cards by {Markov} decision
  processes.
\newblock {\em European Journal of Operational Research}, 212(1):123--130.

\bibitem[Soman and Cheema, 2002]{soman2002creditlimit}
Soman, D. and Cheema, A. (2002).
\newblock The effect of credit on spending decisions: The role of the credit
  limit and credibility.
\newblock {\em Marketing Science}, 21(1):pp. 32--53.

\bibitem[Thomas et~al., 2002]{thomas2002credit}
Thomas, L.~C., Edelman, D.~B., and Crook, J.~N. (2002).
\newblock {\em Credit scoring and its applications}.
\newblock SIAM.

\bibitem[Thomas et~al., 2004]{thomas2004readings}
Thomas, L.~C., Edelman, D.~B., and Crook, J.~N. (2004).
\newblock Readings in credit scoring: foundations, developments, and aims.
\newblock {\em OUP Catalogue}.

\bibitem[Till and Hand, 2003]{till2003behavioural}
Till, R. and Hand, D. (2003).
\newblock Behavioural models of credit card usage.
\newblock {\em Journal of Applied Statistics}, 30(10):1201--1220.

\bibitem[Till, 2001]{tillThesis2001}
Till, R.~J. (2001).
\newblock Predictive behavioural models in credit scoring and retail banking.
\newblock Unpublished PhD Dissertation, Imperial College, London.

\bibitem[Trench et~al., 2003]{trench2003managing}
Trench, M.~S., Pederson, S.~P., Lau, E.~T., Ma, L., Wang, H., and Nair, S.~K.
  (2003).
\newblock Managing credit lines and prices for {Bank One} credit cards.
\newblock {\em Interfaces}, 33(5):4--21.

\end{thebibliography}
\addcontentsline{toc}{section}{References}
\appendix
\section*{Proof of the scaling property of $\hat{\ell}$}
\label{app:scalingProof}
By the assumption that $A'(t) \eqdist \alpha A(t)$,
\begin{align*}
  \hat{\ell}' & = \argmax_{\ell'} \Big\{ \gamma \expt{B_{\ell'}'(t)} - \nu
                \ell' \Big\} \\
              &  =
                \argmax_{\ell'} \Big\{ \gamma \expt{\sup_{0 \leq u \leq t} \{A'(u) : A'(u) \leq
                \ell'\} } - \nu \ell' \Big\} \\
              & = \alpha \argmax_{\ell'} \Bigg\{ \gamma \expt{ \sup_{0 \leq u \leq t}
                \left\{ A(u) : A(u) \leq
                \frac{\ell'}{\alpha} \right\} } - \nu \frac{\ell'}{\alpha},
                \Bigg\}.
\end{align*}
Making the substitution $\ell = \ell'/\alpha$,
\begin{align*}
  \hat{\ell}' & = \alpha \argmax_{\ell} \Bigg\{ \gamma \expt{\sup_{0 \leq u \leq t}
                \left\{ A(u) : A(u) \leq
                \ell \right\} } - \nu \ell \Bigg\} \\
              &= \alpha \argmax_{\ell} \Big\{ \gamma \expt{B_\ell(t)}  - \nu \ell
                \Big\} = \alpha \hat{\ell},
\end{align*}
which shows that the optimal limit scales with $\alpha$. Since
$\hat{\ell}' = \alpha \hat{\ell}$, we have
\begin{align*}
  \prob{A'(t) > \hat{\ell}'} = \prob{\alpha A(t) > \alpha \hat{\ell}} = \prob{A(t) >
  \hat{\ell}},
\end{align*}
which shows that the blocking probabilities remain the same.
\end{document}